\input amstex
\documentstyle{amsppt}
\input psfig.sty
\def\s{\vskip6pt}
\def\codim{{\roman codim}}
\def\mpr#1{\;\smash{\mathop{\hbox to 20pt{\rightarrowfill}}\limits^{#1}}\;}
\def\epi#1{\;\smash{\mathop{\hbox to 20pt{\rightarrowfill}\hskip
-13pt\rightarrow}\limits^{#1\,}}\;}
\def\epii{\smash{\mathop{\hbox to 14pt{\rightarrowfill}\hskip
-11pt\rightarrow}}}
\def\mono{\lhook\joinrel\relbar\joinrel\rightarrow}
\def\mpl#1{\;\smash{\mathop{\hbox to 20pt{\leftarrowfill}}\limits^{#1}}\;}

\def\dea{\downarrow\kern-1.12em\lower 2pt\hbox{ $\downarrow$}}
\def\swea{\swarrow\hskip-17.75pt\lower 1.4pt\hbox{ $\swarrow$}}

\def\I{{\Cal I}}
\def\K{{\Cal K}}
\def\Ct{{\Bbb C}}
\def\Zt{{\Bbb Z}}

\def\D{{\Bbb D}}

\def\CP{{\Bbb{ CP}}}
\def\Hi{H^\flat}
\def\Ho{H^\sharp}
\def\Pii{P^\flat}
\def\Po{P^\sharp}

\def\ind{\widetilde{\kappa}}
\def\Ind{\kappa}

\def\P{{\Cal P}}
\def\F{{\Cal F}}

\topmatter
\title
Generalized Riemann-Hilbert Transmission\\ and Boundary Value
Problems,\\ Fredholm Pairs and Bordisms
\endtitle
\author Bogdan Bojarski and Andrzej Weber \endauthor
\thanks
Both authors are supported by the European Commission RTN
HPRN-CT-1999-00118, {\it Geometric Analysis}.  The second author is
supported by KBN 2P03A 00218 grant. The second author also thanks
Instytut Matematyczny PAN for hospitality. 
\endthanks

\abstract We present classical and generalized Riemann-Hilbert problem in 
several
contexts arising from $K$-theory and bordism theory. The language of
Fredholm pairs turns out to be useful and unavoidable. We propose an
abstract formulation of a notion of bordism in the context of Hilbert
spaces equipped with splittings.
\endabstract                                    
\endtopmatter
\document
\rightheadtext{ Riemann-Hilbert Problem, Fredholm Pairs, Bordisms}
\head \bf\S1.  Introduction \endhead
The concept of a Fredholm pair $\P=(H^-,H^+)$ of closed subspaces
$H^-$, $H^+$ of a Hilbert (or Banach) space was introduced in 1966 by T.
Kato in his studies of stability properties of closed, mainly
unbounded operators, \cite{19}.

Recall that the pair $\P=(H^-,H^+)$ is a Fredholm pair if
the algebraic sum $H^-+H^+$ is closed and the numbers
$\alpha_\P=\dim(H^-\cap H^+)$ and
$\beta_\P=\codim(H^-+H^+)$ are both finite. We also assume that $H^-$ and
$H^+$ are of infinite dimensions. The difference $\alpha_\P-\beta_\P$
was defined in \cite{19} as the index of the pair, $Ind\,\P$, and the
crucial observation of T. Kato was that $Ind(H^-,H^+)$ is not changed by
,,small'' deformations of the pair. More precisely, the set $\F
Gr^2(H)$ of all Fredholm pairs of a Hilbert space appears then as
an open subset of the Cartesian product $Gr(H)\times Gr(H)$ of
the Grassmannian of closed subspaces of  $H$  supplied  with  the
usual
,,minimal gap'' metric, see \cite{17, 19}. In this context the notation $\F
Gr^2(H)$ can be interpreted as the Fredholm bi-Grassmannian of the
Hilbert space $H$ which generalizes in a natural way to the Fredholm
multi-Grassmannian $\F Gr^n(H)$, when instead of pairs of
subspaces we consider $n$-tuples $(M_1,\dots,M_n)$ of closed
subspaces forming Fredholm fans, \cite{4, 6}.

There was no doubt from the outset that the theory of Fredholm pairs
and their generalizations should be studied in close
relationship with the theory of Fredholm operators. Thus
Fredholm pairs in Kato's \cite{19} were considered as a convenient
extension of the theory of Fredholm operators.  For a
Fredholm, possibly unbounded, closed operator $A:H_1\rightarrow
H_2$ acting between Hilbert (Banach) spaces, the pair
$\P_A=(graph\,A,\widetilde H_1)$ was in \cite{19} the basic example
of a Fredholm pair. Here the ,,coordinate'' subspace $\widetilde
H_1=H_1\oplus 0$ and the graph of $A$ are closed subspaces in the direct
sum $H=H_1\oplus H_2$.
Moreover, $$Ind\,\P_A=ind\,A\,$$ where $ind\,A$ denotes here and
in the sequel the index of the Fredholm operator $A$.
Also in \cite{4} the bi-Grassmannian $\F Gr^2(H)$, understood
there also as the space of abstract Riemann-Hilbert transmission
problems, was parameterized by a family of Fredholm operators ${\Cal L}_\P$
associated with projectors $(P^-,P^+)$, not necessarily orthogonal,
onto the spaces of the pair $\P$.

The theory of Fredholm operators in Hilbert space turned out to be
an important tool for studying topology of manifolds and
$K$-theory, especially the geometrical and topological invariants
defined by elliptic differential and pseudodifferential operators
in spaces of sections of smooth vector bundles on manifolds. The
highlight along that road was the famous solution by M.~Atiyah and
I.~Singer, \cite{2}, of the index problem for elliptic operators.
In the abstract functional analytic setting the space $\F (H)$ of
Fredholm operators in the Hilbert space $H$, topologized as a
subset of the Banach algebra ${\Cal B}(H)$ of bounded operators in
$H$, turned out to be the classifying space for the functor
$K^0(-)$, the 0-th term of the generalized cohomology theory
$K^*(-)$, \cite{1}. Later the $K$-homology $K_*(X)$ of a
topological space $X$ (or $K^*(A)$ for a $\Ct^*$-algebra in the
noncommutative case) was introduced, \cite{18}. According to
Kasparov the generators of $K_*(X)$ are realized by certain
Fredholm operators acting in Hilbert space, which is equipped with
an action of the algebra of functions $C(X)$.

The roots of the
extremely successful applications of the Fredholm operators in
global analysis, geometry of elliptic operators and $K$-theory,
undoubtedly are related with the following basic features of the class
$\F (H)$:\widestnumber\item{(iii)}
\s
\item{ (i)} The set $\F (H)$ is stable under sufficiently small perturbations
in ${\Cal B}(H)$ i.e.
$$A\in \F (H)\Rightarrow A+B+K\in \F (H)$$
for $B\in {\Cal B}(H)$, $\|B\|<\varepsilon_A $ (for a sufficiently small
$\varepsilon_A$) and $K\in \K(H)$, where $\K(H)$ denotes the ideal of compact
operators in $H$;
\s\item{ (ii)} Composition law: $$A\in \F (H)\,,\; B\in
\F (H)\Rightarrow A\circ B\in \F (H)$$
and $ind\,A\circ B=ind\,A+ind\,B\,.$ The index homomorphism
$$ind:\F (H)\rightarrow \Zt$$ is surjective and describes the
set of components $\pi_0(\F (H))$;
\s\item{ (iii)} In interesting and important cases, arising in
the theory of partial differential equations and boundary value
problems, the Hilbert space appears as a function space over a
manifold, usually a function space of Sobolev type. Therefore it was
natural to consider
Hilbert spaces equipped with an action of the algebra of functions $B=C(X)$
over a topological space (usually a manifold) $X$.
More generally, it was assumed in \cite{18} that the considered Hilbert spaces
are
Hilbert modules with some $\Ct^*$-algebra $B$ action
$$r:B\rightarrow {\Cal B}(H)\,.$$
The condition
$$\forall b\in B\,:\;[r(b),A]\in\K$$
distinguishes a class of operators $A\in \Cal B$ which is of special
interest.
In consequence it restricts also the class of Fredholm operators.
It is a remarkable fact, that the elliptic pseudodifferential
operators belong to the described above class for the standard
multiplication representation of the algebra of continuous functions.
\s

\noindent
The calculus of commutators and their traces was the starting point
for A.~Con\-nes for introducing cyclic cohomology and
proclaiming the program of noncommutative geometry, \cite{12}.

The natural and intimate connection of the classical Riemann-Hilbert
transmission problems and the theory of Fredholm pairs in a Hilbert
space $H$ was first described in middle seventies by the first
author, \cite{4}.
In particular the mentioned above basic properties (i), (ii) and (iii)
appear in a decisive way in \cite{4}.
In the linear transmission problems for
Cauchy-Riemann systems, generalized Cauchy-Riemann systems, Dirac
operators as well as higher dimensional transmission problems related
with the Cauchy data spaces for higher order elliptic operators in vector
bundles on manifolds, the Fredholm pair approach is more direct then
the usual reduction process to systems of elliptic $\Psi$DO's on the
boundary or the splitting submanifold. In \cite{4} a variety of concepts
have  been  introduced.  Besides   the   named   above   Fredholm
bi-Grassmannian
and the abstract Riemann-Hilbert transmission problem let's mention
here the discussion of the role of Calder\'on projectors on the
Cauchy data spaces in general vector bundle setting, Green
formulas and pairing between Cauchy data spaces for $D$ and the
formally adjoint operator $D^*$, splitting index formulas. In the case
when $H=H_1\oplus H_2$ Fredholm pairs were discussed as pairs of
correspondences (relations), which may be composed, leading to
a generalization of composition rules for Fredholm operators. This is
crucial for the case of bordisms, \cite{8} and \S5 
below.

Applications to topology of Fredholm pairs are not enough investigated
so far. Except for the articles [Bo1-4], there are very
few papers exploring this subject. One should mention
\cite{9-10}.

The purpose of this note and its expanded version \cite{8}
is to give an introduction to a systematic treatment of the Fredholm
pairs theory applied to geometry and topology.
In terms of
boundary values of solutions the Riemann-Hilbert
problem translates directly to the language of Fredholm pairs.
One
can recover the index of the original problem as well as the kernel
and the cokernel. Developing this idea we study an application of
Fredholm pairs to bordisms. We
consider a bordism of smooth manifolds $$M_1\sim_X M_2$$ equipped
with an elliptic differential operator $D$ acting on the sections of a vector
bundle $\xi$
over $X$. Let $H_i=L^2(M_i;\xi)$.
The generalized boundary values of the solutions of $Du=0$
form a subspace $L$ contained in the direct sum $H_1\oplus H_2$.
This space cannot be represented as a graph of an operator
$H_1\rightarrow H_2$, but it may be
treated as a morphism from $H_1$ to $H_2$. It transports certain
family of
linear subspaces of $H_1$ to $H_2$. The correspondence $L$
allows to couple spaces $H_1^-\subset H_1$ with spaces
$H_2^+\subset H_2$. In general the index is defined when
$(L,H_1^-\oplus H_2^+)$ is a Fredholm pair.

As in the case of Fredholm operators, the properties (i)--(iii) play
the decisive role:
\s\item{ (i)} The index of Fredholm pairs is stable under
deformations;
\s\item{ (ii)} Although the index is not additive under the composition
of correspondences, but the defect is well understood;
\s\item{ (iii)} The constructions are motivated and illustrated by
examples coming from the boundary value problems of elliptic
operators.\s

\noindent
It appears that the concept of Fredholm pairs and correspondences
creates a natural analytical setting for an abstract theory of
bordisms, expressed in terms of linear functional analysis.

The concept of Fredholm pairs and its generalizations provide a
convenient approach to a variety of problems in partial
differential equations: both local, as classical boundary value
problems, or non-local, when in the boundary conditions a global
operator (e.g.~spectral restriction or additional pointwise translation) is
present. In the setting of Fredholm pairs the given
differential operators and their parametrices exist on the same
footing. Families of Fredholm pairs and the bi-Grassmannian $\F
Gr^2$ appear as a classifying space for $K$-functor.
The algebraic construction of $K$-homology $K_*(X)$ suggested by Atiyah
and realized by Kasparov, based on the theory of elliptic or
Fredholm operators in $C(X)$ (or $C^\infty(X)$) modules, have direct
analogies in the Fredholm pairs setting. Some construction e.g.
description of the differential
$$\delta:K_0(X)\rightarrow K_1(M)$$ in the $K$-homology for the
Mayer-Vietoris exact sequence for a splitting $X=X_-\cup_M X_+$ is
easier then in the Fredholm operator setting, \cite{8}. In some
situations, like the Cauchy data spaces for elliptic operators or
the bordism category, the language of Fredholm pairs and
correspondences seems unavoidable. The concept of abstract
Fredholm pairs and bordisms admits natural and well motivated
generalizations: Fredholm fans, \cite{6}, also treated in
\cite{8}.

\head\bf \S2. 
Classical and abstract Riemann-Hilbert pro\-blem  \endhead

The classical Riemann-Hilbert problem is understood as follows. Let
$\CP^1=\D_-\cup_{S^1}\D_+$ be the usual decomposition of the Riemann
sphere (i.e.~the complex projective line). Here $\D_+$ is the unit
disk and $\D_-$ is the complementary disk containing infinity.
Given a function (a loop) $\phi:S^1\rightarrow GL(\Ct^n)$, describe
the totality of
holomorphic vector-valued functions $s_\pm:\D_\pm\rightarrow\Ct^n$,
such that $s_+(z)=\phi(z)s_-(z)$ for $z\in S^1$.
Due to Birkhoff decomposition, \cite{22, 26}, if $\phi$ is
differentiable then the Riemann-Hilbert problem is the same as
looking for a section of the holomorphic bundle defined by $\phi$,
\cite{21, 5, 7, 9, 10}. If
$\phi$ is piecewise constant, then this is the $21^{\roman{ st}}$  problem, as
stated by Hilbert, see \cite{11}.  It's a question about
existence of a system of
singular differential equations with prescribed monodromy.

Denote by $H^\pm$ the space of boundary values of holomorphic
vector-functions on $\D_\pm$. This is a Fredholm pair in
$H=L^2(S^1;\Ct^n)$. The pair
$(\phi(H^-),H^+)$ is also Fredholm, where $\phi(H^-)$ is the image of $H^-$
with respect to the obvious multiplication representation of the loop group
$\Omega
GL(\Ct^n)$ on $H$.

If we normalize $s_-$ by the condition
$s_-(\infty)=0$, then we obtain a subspace $\Hi\subset H^-$. Set
$\Ho=H^+$. Then $H=\Hi\oplus\Ho$.
According to \cite{4}
the question about the pair $(\phi(H^-),H^+)$ being Fredholm
is reduced to the abstract problem of studying the operator
$${\Cal L}_\phi=\phi\Pii+\Po:H\rightarrow H$$
or the Toeplitz operator
$$\Pii\phi:\Hi\rightarrow\Hi\,.$$
The projectors $\Pii$ and $\Po$ are the projectors in the direct sum
$H=\Hi\oplus\Ho$, but they can be substituted by the
projectors of Sohotski-Plemelj, \cite{25, 15}, which are singular integral
operators.

In order to explain the deep meaning of the operator ${\Cal L}_\phi$
let us summarize some facts about Fredholm pairs. We will
follow \cite{4} and \cite{6}.
Suppose that $H$ is decomposed into a direct sum
$$H=\Hi\oplus \Ho\,,$$ both summands closed of infinite dimension.
We can say that this decomposition is given by a symmetry $S$:~a
,,sign'' or ,,signature'' operator with $S^2=1$. Let $\Pii$ and
$\Po$ be the corresponding projectors and $S=\Po-\Pii$. This is
the basic one-dimensional singular integral operator. It is well
known that for any continuous loop $\phi$ the commutator
$[\phi,S]=\phi S-S\phi$ is a compact operator (Mikhlin lemma,
\cite{26}).

Let ${\Cal I}\subset {\Cal B}(H)$ be an ideal containing the ideal
of finite rank operators and contained in the ideal of compact
operators.  Define $GL(S,{\Cal I})\subset GL(H)$ to be the set of all
invertible authomorphisms of $H$ commuting with $S$ up to the ideal
$\Cal I$: $$GL(S,\I)=\{\phi\in GL(H)\,:\;[\phi,S]\in\I\}\,.$$
We have the
following classification result.

\proclaim {Theorem 2.1}{\rm\cite{4}}
Let $H^\pm$ be a Fredholm pair
with $H^+=\Ho$.
Suppose, that $H^-$ is given by a projector $P^-$ satisfying
$\Pii-P^-\in \I$.
Then there exists $\phi\in GL(S,{\Cal I})$, such that
$H^-=\phi(\Hi)$.
Moreover, the operator ${\Cal L}_\phi=\phi\Pii+\Po$ is Fredholm and
$$ind({\Cal L}_\phi)=Ind(H^-,H^+)\,.$$
The map
$$\ind :GL(S,{\Cal I})\rightarrow \Zt$$
$$\ind (\phi)=ind({\Cal L}_\phi)$$
is a group homomorphism
$$\ind(\phi\circ\psi)=\ind(\phi)+\ind(\psi)\,.$$
\endproclaim

We remark that the correspondence $\phi \mapsto {\Cal L}_\phi$
between the group $GL(S,\I)$ and the Fredholm operators is an ,,almost''
homomorphism,  i.e.
$${\Cal L}_{\phi\circ\psi}={\Cal L}_{\phi}\circ
{\Cal L}_{\psi}+T(\phi,\psi)$$
with $T(\phi,\psi)={1\over  4}(1-\phi)[S,\psi](1-S)\in \I$.
If we write $\phi\in GL(S,{\Cal I})$ in a matrix form with respect to
the splitting $H=\Hi\oplus\Ho$:
$$\phi=\left(\matrix\alpha&\beta\\ \gamma& \delta \endmatrix\right)\,,$$
then $\alpha$, $\delta$ are Fredholm operators, $\beta$, $\gamma$ are in $\I$
and $\ind (\phi)=ind(\alpha)=-ind(\delta)$.
It follows that
$$Ind(H^-, H^+ ) = ind(\Pii\phi:\Hi\rightarrow
\Hi)=ind(\Po\phi^{-1}:\Ho\rightarrow \Ho)\,.$$
Note, that if $\phi$ is arbitrary, possibly not invertible, then
$\phi\Pii+\Po:H\rightarrow H$ and $\Pii\phi:\Hi\rightarrow\Hi$ are
Fredholm operators of equal indices, provided that
$\alpha$ is Fredholm. These indices are not necessarily equal to
$Ind(\phi(\Hi),\Ho)$. The equality holds if and only if $\phi_{|\Hi}$
is injective. It is
better to distinguish between the domain of $\phi$ (we write $H_1$)
and its target ($H_2$). There is another expression for
$\ind (\phi)$, which will be useful later:
\proclaim{Proposition 2.2}
The pair
$(graph\,\phi,\Hi_1\oplus\Ho_2)$ in $H_1\oplus H_2$ is Fredholm and
its index is equal to $\ind (\phi)$.\endproclaim

We have introduced a splitting of the Hilbert space
$H=L^2(S^1;\Ct^n)=\Hi\oplus\Ho$. It's a good moment now to expose its
fundamental role. The splitting comes from the
division of the Fourier exponents into subsets
$$\Zt=\Zt_{<0}\cup\Zt_{\geq 0}\,.$$
A finite perturbation of this set is also an admissible
decomposition.
The need of introducing a splitting was clear already in
\cite{4}:
\s\item{$\bullet$} It was used to the study of Fredholm pairs with
application to Riemann-Hilbert problem in \cite{4}.
\s\item{$\bullet$} Splitting also came into light
in the paper of Kasparov, \cite{18} who introduced homological $K$-theory
built from Hilbert modules. The program of noncommutative
geometry of A.Connes develops this idea.
\s\item{$\bullet$} Splitting plays an important role in the theory of loop
groups in \cite{26}.
\s\item{$\bullet$} There is also a number of papers in which surgery of the
Dirac operator is studied. Splitting serves as a boundary
condition, see e.g.~\cite{14, 27}. These papers
originate from \cite{3}.
\s

Let us come back to the decomposition of $L^2(S^1;\Ct^n)$ originating from the
classical Riemann-Hilbert problem. It is given by a pair of
pseudodifferential projectors.
Suppose that the authomorphism $\phi$ is the multiplication by a matrix with
entries being continuous functions. Then $\phi\in GL(S,\K)$ and the
Theorem 2.1. 
applies.

\head \bf\S3. Fredholm bi-Grassmannian  \endhead

We will describe the homotopy types of the spaces involved in our
constructions.

\s
\noindent\bf 3.1. \it The
Grassmannian of the closed linear subspaces $M\subset H$ with
$\dim (M)=\codim(M)=\infty$. \rm We denote this set by
$Gr_\infty(H)$. The linear group $GL(H)$ acts on it transitively.
Let $S\in GL(H)$ be a symmetry decomposing $H$ into direct sum
$\Hi\oplus\Ho$ of closed subspaces of infinite dimensions.
The stabilizer of $\Ho\in
Gr_\infty(H)$ consists of linear isomorphisms, which can be
written in the block form $\left(\matrix \alpha&0\\ \gamma&\delta
\endmatrix \right)$
with $\alpha$, $\delta$ being isomorphisms and $\gamma$ arbitrary
linear map. We can write $Gr_\infty(H)=GL(H)/Stab(\Ho)$.
We endow this set with the quotient topology. By a result of
Kuiper, \cite{20}, the topological spaces $GL(H)$ and $Stab(\Ho)=GL(\Hi)\times
Hom(\Hi,\Ho)\times GL(\Ho)$ are contractible. Hence $Gr_\infty(H)$
is contractible as well.

\s
\noindent\bf 3.2. \it The set of Fredholm pairs $(H^-,H^+)$ in $H$. \rm
We denote this set by
$\F Gr^2(H)$. This is a subset of $Gr_\infty(H)\times Gr_\infty(H)$.
The projection on the second factor (forgetting about $H^-$) is a
fibration. Denote the fiber over $\Ho$ by $Gr_{\Ho}(H)$
$$ Gr_{\Ho}(H)\mono \F Gr^2(H) \epii Gr_\infty (H)\,.$$
Since
the base of the fibration is contractible by \S3.1, 
the inclusion of the fiber is a homotopy equivalence.
\s

\noindent\bf 3.3. \rm The fiber $Gr_{\Ho}(H)$ is identified with the subset of
$Gr_\infty (H)$ consisting of \it the closed linear subspaces $H^-\subset H$,
such that the pair $(H^-,\Ho)$ is Fredholm. \rm It is the orbit of
$\Hi$ with respect to the action of the ,,parabolic up
to $\K$'' subgroup $P(S,\K)\subset GL(H)$
$$Gr_{\Hi}(H)=P(S,\K)\cdot\Hi\subset Gr_\infty(H)\,.$$
The group $P(S,\K)$ consists of isomorphism
of the form $\left(\matrix \alpha&\beta \\ \gamma&\delta \endmatrix\right)$ 
with
$\alpha$, $\delta$ being Fredholm operators and $\beta$ compact
operator. Consider the projection $P(S,\K)\rightarrow \F (\Hi)$
sending an element of $\phi\in P(S,\K)$ to the operator
$\alpha=\Pii\phi_{|\Hi}$.
Arguing as in \cite{26} 6.2.4 we prove that this map has the contractible
fibers. Therefore it is a homotopy equivalence.

Let us sum up the results and compare it with \cite{26}, where the
restricted Grassmannian was studied.

\proclaim {Theorem 3.4} The following maps are homotopy equivalences:
$$\matrix
GL(S,\K )&\mono&P(S,\K)&\epii&\F (\Hi)\\
\dea & & \dea & &\phantom{\frac \int\int}\\
 Gr_{res}(S,\K)&\mono&  Gr_{\Ho} (H)&\mono& \F Gr^2(H)\,,\\
\endmatrix$$
where $ Gr_{res}(S,\K)$ consists of the subspaces $H^-\subset H$ of
the form $H^-=\phi(\Hi)$ for $\phi\in GL(S,\K)$.\endproclaim

In our constructions the ideal of compact operators $\K$ can be
substituted by any smaller ideal containing the ideal of the finite
rank operators. In fact in \cite{26} the ideal of Hilbert-Schmidt
operators appears.

\head\bf \S4. Some general remarks and comments  \endhead

We recapitulate: there are several aspects of the general
Riemann-Hilbert problem.
\s

\noindent\bf 4.1. \it The underlying geometric and combinatorial setting. \rm 
The
simplest example of the Riemann-Hilbert problem deals with
the operator $\overline\partial$ on the Riemannian sphere
$S^2=\CP^1$ (or a Riemann surface) divided into two complementary
domains by a curve $\Gamma$, say $\Gamma=S^1$, as described. In
the literature much more general situations have been studied,
\cite{25, 15, 22}.

Suppose that there is given an oriented contour $M$ on a
Riemann surface $X$. It may consist of a finite number of
parameterized curves admitting transversal intersections and
self-intersections. We have on each component of $M$ positive and
negative side (i.e.~an orientation of the normal bundle).

\centerline{\psfig{file=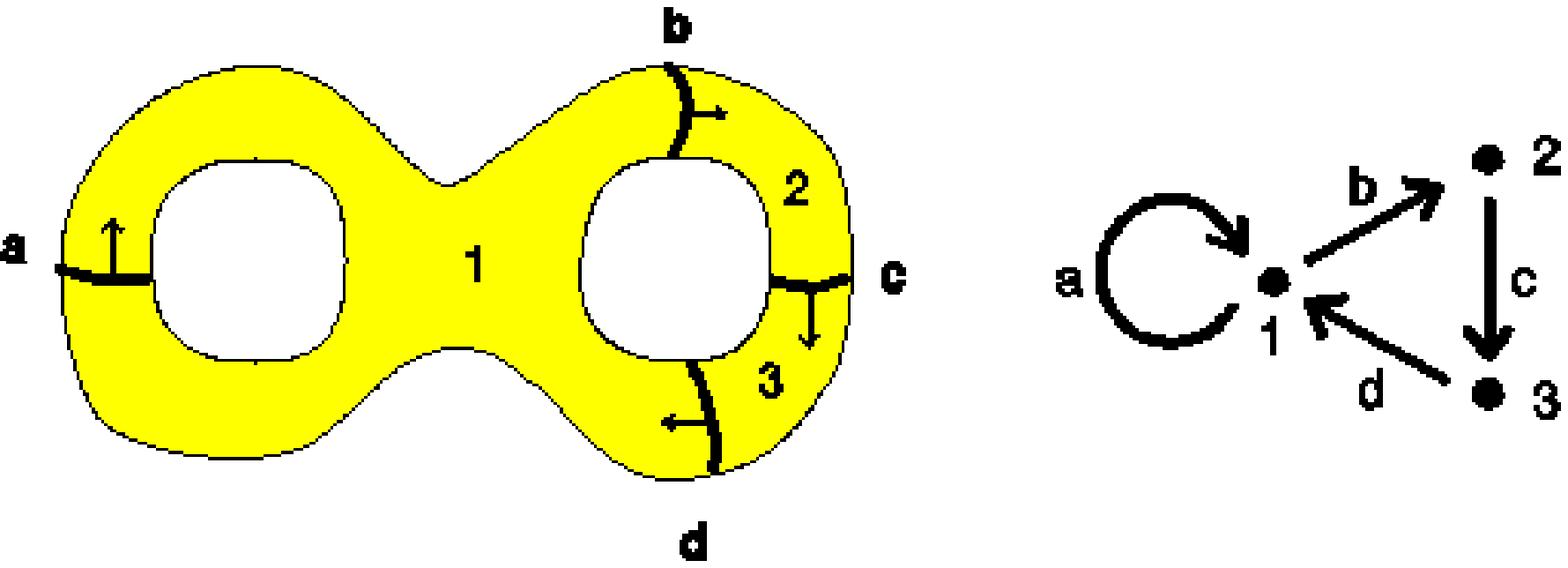,height=1truein}}

\centerline{\tt Configuration of curves and its combinatorial
model.} \noindent
If $$\phi:M\rightarrow GL(\Ct^n)$$
is a continuous map, then for a vector-function $s(z)$ holomorphic
on each component $U$ of $X\setminus M$, admitting in some natural
sense boundary extension to $\overline U$, the Riemann-Hilbert
transmission condition $$s_+(z)=\phi(z)s_-(z)\,,\quad{\roman{ for}}\quad z\in 
M$$
 is meaningful. Here $s_\pm(z)$ are boundary values of $s(z)$ on positive and
negative side of $M$.

If $M=\bigcup_{i=0}^m \Gamma_i$ is the sum of boundaries of
disjoint discs $\Gamma_i=\partial \D_i$, $i=0$,\dots $m$, then we
have the classical Riemann-Hilbert problem in a non-simplyconnected
domain. Splitting $M$ into two disjoint parts $M=M_1\sqcup M_2$ provides
an example of a bordism. Such bordisms admit compositions, if the
boundaries are matching. One can formulate local and global index
formulas in the realm of conformal field theory. For details see \S6
and \cite{8}.

In higher dimensions one considers manifolds with a
configurations of hypersurfaces (of codimension one)
non-intersecting or with transversal intersections. The most
relevant here is the case
of submanifolds realizing a decomposition of $X$ into bordisms
$$X_0\cup_{M_1}X_1\cup_{M_2}\dots\cup_{M_m}X_m\,.$$

Recall that, as beautifully described in \cite{24}, the bordisms
form a category with oriented $n-1$-dimensional manifolds as objects
and bordisms as morphisms. In our abstract bordism model each splitting
manifold $M_i$ has an associated Hilbert space $H_i$ supplied with an
involution $S_i$ also called signature operator. These define
splittings $H_i=\Hi_i\oplus\Ho_i$ into incoming and outgoing
components and should be considered as a part of the structure. The
bordism $M_{i-1}\sim_{X_1}M_i$ together with an elliptic first order operator
$D$ on $X_i$ gives rise to a closed linear subspace
$L\subset H_{i-1}\oplus H_i$,
the space of Cauchy data on $\partial X_i=M_{i-1}\sqcup M_i$ of
solutions of homogeneous equation $D=0$ on $X_i$.
The space $L$ is a correspondence from $H_{i-1}$ to $H_i$.
We will illustrate our point of view by a simple but instructive
Example 5.1.
\s

\noindent\bf 4.2. \it The generalized Riemann-Hilbert problem. \rm
In the simplest case $m=1$
$$X_-\cup_{M}X_+$$
Fredholm pair arises from consideration of the Cauchy data on
$X_-$ and $X_+$. More precisely let $$ D
:C^\infty(X;\xi )\rightarrow C^\infty(X;\eta)$$ be an elliptic
operator of the first order.
The Dirac operator is of special interest.
One defines the spaces $H^\pm( D
)\subset H=L^2(M;\xi)$, which are
the spaces of boundary values of solutions of homogeneous equations $Ds=0$ on
the
manifolds $X_\pm$. The pair $H^\pm(D)$ is Fredholm.
In order to study $ker\,D$ and $coker\,D$ it
is convenient to assume that $D$ and $D^*$ have
the unique extension property, i.e.~the solutions of $D$ and $D^*$
are determined by the boundary values on $M$. Then
$$ker\,D= H^+(D)\cap H^-(D)\,,\quad coker\,D= H/(H^+(D)+ H^-(D))\,.$$
As in the
case described in \S2 
corresponding Cauchy data spaces
admit projectors. There are associated Calder\'on projectors
$P^\pm(D)$ onto $H^\pm(D)$. They are complementary up to a compact
operator: $P^-(D)+P^+(D)-1\in \K\,.$ The operation $S$ given by
$$S=P^+-P^-$$ is the fundamental singular operator.
The group $GL(S,\K )$ is naturally involved.
\s

\noindent\bf 4.3. \it The Riemann-Hilbert problem in an abstract Hilbert space
$H$: \rm  suppose we have an involution $S\in{\Cal B}(H)$ defining a splitting
$\Hi\oplus\Ho$. We consider the following objects:

\s\item{$\bullet$} the group of $GL(S,\K )$ of the linear isomorphisms of $H$
commuting with $S$ up to $\K$,

\s\item{$\bullet$} the bi-Grassmannian $\F Gr^2_{res}$ consisting of pairs of 
the form
$(\phi(\Hi),\Ho)$, with $\phi\in GL(S,\K )$, called restricted
Grassmannian  in \cite{26}.

\s\item{$\bullet$} the bi-Grassmannian $\F Gr^2(H)$, the set of all Fredholm 
pairs in $H$.
\s
\noindent
These spaces are homotopy equivalent, they are classifying
spaces of $K$-theory $BU\times\Zt$.
A family of Fredholm
pairs (say over $T$) defines an element of $K^0(T)$.
Moreover, the classical Riemann-Hilbert problem gives us a way
of constructing a Fredholm pair in $H^n=L^2(S^1;\Ct^n)$ out of
a given loop in $U_n\subset GL(\Ct^n)$. The assignment
$$\Omega U_n\rightarrow \F Gr^2(H^n)\simeq BU\times\Zt$$ passes to a
map
$$\Omega U_\infty\rightarrow \F Gr^2(H^\infty)\simeq BU\times\Zt\,,$$
which can be interpreted as the Bott periodicity map.
\s

\noindent\bf 4.4. \it Quantum Riemann-Hilbert problem: \rm
There is another structure which one cannot forget keeping in mind
geometric applications. The Hilbert space $H$ comes with an action of
the algebra $C(M)$ of functions on $M=\partial X_\pm$. The
pseudodifferential operator $P^+(D)-P^-(D)$ almost
commutes with the algebra action. Hence it defines an element in the
odd $K$-homology $K_1(M)$.
From the point of view of Kasparov
theory we can replace $P^+(D)-P^-(D)$ by the almost equal operator
$S=\Po-\Pii$, $S^2=1$. But
now as in \cite{13}, pp 287-289, it is easier to express the
pairing with $K^1(M)$. If $\phi\in GL(S^{\oplus n},\K)$ is defined by
a matrix of functions $\widetilde\phi:M\rightarrow GL(\Ct^n)$
(i.e. $\widetilde\phi$ is a generator of $K^1(M)$), then the effect
of the pairing of $[\widetilde\phi]\in K^1(M)$ with $[S]\in K_1(M)$  is
equal to $\ind(\phi)$. One can ask what element of $K_1(M)$ is
defined by $S$. It's not hard to guess that:

\proclaim{Theorem 4.5} The Fredholm module $(H=L^2(M;\xi),S)$
is the image
of $[D]$ with respect to the differential $\delta: K_0(X)\rightarrow
K_1(M)$ in homological Mayer-Vietoris sequence of the triple
$(X,X_-,X_+)$.\endproclaim

\noindent Since our proof in \cite{8} is obtained by means of duality, the
result holds modulo torsion in $K_1(M)$.

It is clear that the algebra of
functions $C(S^1)$ (or $C(M)$) may be replaced by an arbitrary
$\Ct^*$-algebra, possibly noncommutative. The framework of
noncommutative geometry, \cite{12, 13}, is another possible setup
for studying Riemann-Hilbert problem and corresponding Fredholm pairs.

\head\bf \S5.
Bordisms  \endhead
Now we would like to describe more general objects than the operators
$\phi\in GL(S,\K)$ considered so far. We study relations in $H$ or
correspondences from $H_1$ to $H_2$. Our approach is motivated
by the geometric
theory of bordisms, \cite{24}. First, let us present an example:

\remark{Example 5.1}
\rm For $0<r<R$ consider the ring
$$X=\{z\in\Ct\,:\,r\leq |z|\leq R\}\,.$$
Then $$\partial X=M_1\sqcup M_2= S^1_R\sqcup S^1_r\,.$$
The functions $e_i=z^i:M_1\rightarrow \Ct$ and $\epsilon_i=z^i:M_2\rightarrow
\Ct$ for $i\in\Zt$ form a basis of the Hilbert spaces
\s
\hfil$H_1=L^2(M_1;\Ct)=\left\{\sum_{i\in\Zt}a_ie_i\,:
\,\sum_{i\in\Zt}a_i^2R^{2i}<\infty\right\}$
\s
\noindent and
\s
\hfil$H_2=L^2(M_2;\Ct)=\left\{\sum_{i\in\Zt}b_i\epsilon_i\,:
\,\sum_{i\in\Zt}b_i^2r^{2i}<\infty\right\}\,.$
\s
\noindent Consider the Cauchy-Riemann operator acting
on the complex-valued functions on $X$.
The space of the boundary values of solutions $L$ is the graph of the
unbounded operator
\s
\hfil$\def\awr{\mathrel{\smash -}}
\def\przearrow{\awr\awr\rightarrow}
\Phi :H_1\przearrow H_2\,,\qquad
\Phi\left(\sum_{i\in\Zt}a_ie_i\right)=\sum_{i\in\Zt}a_i\epsilon_i\,.$
\s
\noindent The maximal domain of $\Phi$ is
\s
\hfil $\left\{\sum_{i\in\Zt}a_ie_i\,:
\,\sum_{i\in\Zt}a_i^2(R^{2i}+r^{2i})<\infty\right\}\,.$
\s
\noindent The above condition can be substituted by
$\sum_{i<0}a_i^2r^{2i}+\sum_{i\geq 0}a_i^2R^{2i}<\infty$.
Set
$$\Hi_1=span\{e_i\,:\,i<0\}\,,\quad\Ho_1=span\{e_i\,:\,i\geq0\}\,,$$
$$\Hi_2=span\{\epsilon_i\,:\,i<0\}\,,\quad
\Ho_2=span\{\epsilon_i\,:\,i\geq0\}\,.$$
Now $\Phi$ restricted to $\Ho_1$ is bounded and moreover, it is {\it
compact}. Indeed,
$\Phi_{|\Ho_1}$ is given by Cauchy integral
$$\Phi(f)(\zeta)={1\over 2\pi i}\int_{|z|=R}{f(\zeta)d\zeta\over
z-\zeta}\,.$$
Similarly $\Phi^{-1}$ restricted to $\Hi_2$ is a compact operator.
The space of the boundary values of the solutions is the direct
sum of the graphs
$$L=graph(\phi_1)\oplus graph(\phi_2)\,,$$
where
$$\phi_1=\Phi_{|\Ho_1}:\Ho_1\rightarrow \Ho_2\,,\quad\phi_2=
\Phi^{-1}_{|\Hi_2}:\Hi_2\rightarrow \Hi_1\,.$$
Note, that

\s\item{$\bullet$} The pair $(L,\Hi_1\oplus\Ho_2)$ spans $H_1\oplus H_2$ as a
direct sum,
\s\item{$\bullet$} the projection of $L$ onto $\Ho_1\oplus \Hi_2$ along 
$\Hi_1\oplus
\Ho_2$ is an isomorphism,
\s\item{$\bullet$} $L$ is a direct sum of graphs of the compact operators
$\phi_1=\Phi_{|\Ho_1}$ and $\phi_2=\Phi^{-1}_{|\Hi_2}$.
\s

\noindent
Now consider the Cauchy-Riemann operator $D$ on the projective line $\CP^1$
decomposed into the subsets

$$\aligned
X_1=&\{z\in\Ct\,:\,|z|\geq R\}\cup \{\infty\},\\
X=&\{z\in\Ct\,:\,r\leq |z|\leq R\},\\
X_2=&\{z\in\Ct\,:\,|z|\leq r\}\,.\endaligned$$
We have the spaces of boundary values of holomorphic functions on $X_1$ and
$X_2$
$$\aligned H^-_1=&span\{e_i\,:\,i\leq 0\}=\Hi_1\oplus \langle e_0\rangle,\\
H^+_2=&span\{\epsilon_i\,:\,i\geq 0\}=\Ho_2.\endaligned$$
The space $L\subset H_1\oplus H_2$ is not a
graph of a bounded operator, but as in the case of Riemann-Hilbert
transmission problem, we can write
$$ind(D)=Ind(L(H^-_1),H^+_2)\,,$$ where $L(H^-_1)=\{y\in
H_2\,:\,\exists x\in H_1\,,\, (x,y)\in L\}$.
\endremark

The situation described in the example is quite
general.
Consider a manifold $X$ with boundary, which is the sum of two
components $\partial X=M_1\sqcup M_2$. Let
$D:C^\infty(X;\xi)\rightarrow C^\infty(X;\eta)$ be an elliptic
operator of the first order.
Set $H_i=L^2(M_i;\xi)$ for $i=1,2$ and let $L$ be the closure in
$L^2(\partial X;\xi)=H_1\oplus H_2$ of the space $\{ u_{|\partial
X}\,:\, Du=0\,,\;u\in C^\infty(X;\xi)\}$.
Let $P_L$ be the Calder\'on projector
$$P_L:H_1\oplus H_2 \epii L\,.$$
Let $\widetilde\xi_i$ be the pull back of $\xi$ to
$T^*M_i\setminus\{0\}$.
The symbol $\sigma(P_L)_{|M_i}$ is an endomorphism of the bundle
$\widetilde\xi_i$. Let us choose pseudodifferential projectors
$P_i$ acting on $H_i$ with $\sigma(P_i)=\sigma(P_L)_{|M_i}$. Then
$$P_1\oplus P_2\sim P_L: H_1\oplus H_2\rightarrow H_1\oplus H_2\,.$$
Set $\Po_1=P_1$, $\Pii_2=P_2$. These operators define split
Hilbert spaces $H_i=\Hi_i\oplus \Ho_i$ ($i=1,\,2$). It follows that

\s\item{$\bullet$} the pair $(L,\Hi_1\oplus\Ho_2)$ is a Fredholm pair,
\s\item{$\bullet$} the projection $\Po_1\oplus \Pii_2$ from $L$ onto 
$\Ho_1\oplus \Hi_2$
along $\Hi_1\oplus\Ho_2$ is a Fredholm operator.
\s

\noindent
It can be shown that, as in Example 5.1, 
there are compact
operators:
\s\item{$\bullet$}
$\phi_1$ transforming the restrictions on $M_1$ of some
solutions $Du=0$ to the restrictions on $M_2$
\s\item{$\bullet$} and $\phi_2$ acting in the opposite direction,
\s

\noindent
such that up to finite dimensional perturbation
the space $L$ is equal to the sum of their graphs.

Note, that Calder\'on projectors are well defined if we restrict our
considerations
to the space of smooth sections. That means that we work on a
pre-Hilbert level. To obtain the abstract Hilbert space model the
completion operation should be applied. In the completion process
different Sobolev type
metrics have to be used according to the Sobolev trace type
imbedding theorems. It requires caution and involves some
additional technicalities, which have been skipped here. For details
see \cite{8}. The remarks
above can be clearly seen in our basic cobordism example of the
Cauchy-Riemann operator in the ring domain.

We want to define an index of $D$, regardless of all the
possible choices of
manifolds closing $X$. It will be defined with respect to the
splittings.
The index $Ind(L(\Hi_1),\Ho_2)$ is not stable under a compact
perturbation. If we twist $L$ with an authomorphism of the form
$1+K$, $K\in\K$ the index may change. Instead
it is wiser to consider the pair $(L,\Hi_1\oplus\Ho_2)$.
Its index is stable under such twists. It is worth to say when the
considered indices are equal:

\proclaim{Proposition 5.2} $Ind(L,\Hi_1\oplus\Ho_2)=Ind(L(\Hi_1),\Ho_2)$
provided that both following conditions hold
\item{$\bullet$} $L$ is injective on $\Hi_1$, i.e.~if $(x,y)\in L$ and
$(x',y)\in L$, $x,\,x'\in\Hi_1$, then $x=x'$,
\item{$\bullet$} $\Hi_1+ dom\,L=H_1$, where $ dom\,L=\{x\in 
H_1:\,\exists y\in H_2,\,(x,y)\in L\}$.\endproclaim


Consider again the case of a bordism $X$, this time closed from both sides
by manifolds $X_1$ and $X_2$. That is: there is a closed manifold $Y$
with a first order elliptic operator $D$ and $Y$ is decomposed
$$Y=X_1\cup_{M_1}X\cup_{M_2}X_2\,.$$ By Theorem 2.1 there exist authomorphisms
$\phi_i$ of $H_i$ almost commuting with Calder\'on projectors, such that
$$\Hi_1=\phi_1H^-_1(D)\quad {\roman{ and}}\quad H^+_2(D)=\phi_2(\Ho_2)\,.$$
Then (provided that $D$ and $D^*$ have the unique extension property)
$$Ind(D)=Ind(L,H^-_1\oplus H^+_2)=
\ind(\phi_1)+Ind(L,\Hi_1\oplus\Ho_2)+\ind(\phi_2)\,.$$

We see that $L$ plays a role of the twist $\phi_i:H_i\rightarrow H_i$, but
here $L$ allows us to couple ,,the lower half'' of $H_1$ with ,,the upper
half'' of $H_2$. We can treat it as a morphism\footnote{
A different approach to bordisms,
 based on quantum field theory
point of view, is presented in \cite{28}, Lecture 2. Dirac
operators are of special interest.} from $H_1$ to $H_2$.
Note, that $H_1$ and $H_2$ are not canonically identified. Indeed the
manifolds $M_1$ and $M_2$ joined by the bordism $X$ can be quite
different. There are also two different algebras $C(M_1)$ and
$C(M_2)$ acting. The actions commute with the splittings up to
compact operators. The object described here is an abstract
substitute of a geometric bordism.

\proclaim{Definition 5.3} A {\bf restricted bordism} 
$H_1\overset{L}\to{\rightsquigarrow} H_2$ between
split Hilbert spaces
$H_i=\Hi_i\oplus\Ho_i$ (for $i=1,2$) is a closed linear subspace
$L\subset H_1\oplus H_2$, which is the image of a projector $P_L\sim
\Po_1\oplus\Pii_2$.\endproclaim

The widest class of linear correspondences, which allows us to
define the index is the following:

\proclaim{Definition 5.4}
A {\bf Fredholm bordism} $H_1\overset{L}\to{\rightsquigarrow}
H_2$ between split Hilbert spaces
is a closed linear subspace
$L\subset H_1\oplus H_2$, such that the pair $(L,\Hi_1\oplus\Ho_2)$
in $H_1\oplus H_2$ is Fredholm. The index of $L$ is the index of this
pair. It is denoted by $\Ind (L)$ or $\Ind (\Hi_1|L|\Ho_2)$ to expose the
role of splittings.\endproclaim

Note, that by Proposition 2.2 
the graph of an isomorphism
$\phi\in GL(S,\K )$ is a Fredholm bordism and
$\Ind (graph\,\phi)=\ind (\phi)$.

We can say that the class of Hilbert spaces with involutions
(splittings) and Fredholm bordisms $H_1\sim_L H_2$ form a category,
which may be considered as an abstract functional theoretic
counterpart of of the category of geometric bordisms. Each elliptic
differential operator on any geometric bordism, the
Calder\'on projectors and the corresponding involutions gives rise to a
Fredholm bordism.

\head\bf \S6. Riemann surfaces with boundary  \endhead
Let us consider another example which is classical, now
also studied under the name of conformal field theory.
We consider the Hilbert space of complex functions on the circle:
$H=L^2(S^1;\Ct)$.
Let $Y_g$ be a Riemann surface of genus $g$. Suppose we have $k+l$
disjoined holomorphic disks $\D_i$ ($i=1,\dots,k$), $\D'_j$ ($j=1,\dots,l$)
contained in $Y_g$. Let
$X$ be the complement of the disks. We think of $X$ as of a
bordism between $k$ circles and $l$ circles. Let $L\subset H^k\oplus
H^l$ be the space of boundary values of the holomorphic functions on $X$.
Denote by $H(\lambda)$ (for $\lambda\in\Zt$) the space $H$ equipped
with the splitting
$$\Hi(\lambda)=z^\lambda\Hi=span\langle z^i\,:\,i<
\lambda\rangle\,,$$ $$
\Ho(\lambda)=z^\lambda\Ho=span\langle z^i\,:\,i\geq
\lambda\rangle\,.$$
For sequences of integers $\lambda_\bullet=(\lambda_1,\dots,\lambda_k)$ and
$\mu_\bullet=(\mu_1,\dots,\mu_l)$ we have splittings
$$H_1=H^k(\lambda_\bullet)=
\Bigl(\Hi(\lambda_1)\oplus\dots\oplus \Hi(\lambda_k)\Bigr)\oplus
\Bigl(\Ho(\lambda_1)\oplus\dots\oplus \Ho(\lambda_k)\Bigr)
\,,$$ $$
H_2=H^l(\mu_\bullet)=
\Bigl(\Hi(\mu_1)\oplus\dots\oplus \Hi(\mu_k)\Bigr)\oplus
\Bigl(\Ho(\mu_1)\oplus\dots\oplus \Ho(\mu_k)\Bigr)
\,.$$
We will compute the index of $L$ with respect to these
splittings.
An element of the intersection
$L\cap (\Hi_1\oplus\Ho_2)$ defines a meromorphic function on $Y_g$
with zeros (resp.~poles) at the centers of $\D_i$'s (resp. $\D'_j$'s)
of the order at least $\lambda_i$ (resp.~smaller then $\mu_j$).
This is a section of a sheaf
\s
\hfil${\Cal O}\left(-\sum_{i=1}^k\lambda_id_i+
\sum_{j=1}^l(\mu_j-1)d'_j\right)\,.$
\s
\noindent Here $d_i$ and $d'_j$ are the
centers of the disks. The index is equal to the Euler characteristic
of $Y$ with coefficients in this sheaf, that is
\s\hfil
$1-g-\sum_{i=1}^k\lambda_i+\sum_{j=1}^l\mu_j-l\,.$\s\noindent
In particular, if we want to compute the index of the Cauchy-Riemann
operator on $Y$, then we set $\lambda_i=0$, $\mu_j=1$.  These
splittings  agree with the spaces of boundary values of
solutions on the disks: the index is $1-g$. Again we can write 
$Ind(L,\Hi_1\oplus
\Ho_2)=Ind(L(\Hi_1),\Ho_2)$. This number is denoted by $\Ind (\Hi_1|L|\Ho_2)$
according to Definition 5.4.

Note, that the bordisms, as well as correspondences can be composed.
We consider the composition with the splittings coinciding.
If we deal only with connected surfaces then $$\Ind (L_1\circ
L_2)=\Ind (L_1)+\Ind (L_2)\,.$$ If we admit disconnected bordism, then it
may happen, that a closed component is created while sewing the
bordisms. The defect $\Delta=\Ind (L_1)+\Ind (L_2)-\Ind (L_1\circ L_2)$
equals to the index on this component. This remark generalizes to an
arbitrary elliptic differential operator $D$ of the first order. In
consequence a decomposition of a closed manifold
$$X=\emptyset\sim_{X_0}M_1\sim_{X_1}\dots\sim_{X_{n-1}}M_n\sim_{X_n}
\emptyset$$
gives rise to a sequence of restricted bordisms
$$0\overset{L_0}\to{\rightsquigarrow}H_1\overset{L_1}\to{\rightsquigarrow}
\dots\overset{L_{n-1}}\to{\rightsquigarrow}H_n
\overset{L_n}\to{\rightsquigarrow}0\,.$$

\proclaim{Theorem 6.1} {\rm \cite{8}} Suppose $D$ and $D^*$ have the unique
extension property. Fix splittings $S_i$ of $H_i$. Then the global
index of $D$ is equal to the sum of partial indices:
$$ind\,
D=\sum_{i=0}^{n}\Ind(\Hi_i|L_i|\Ho_{i+1})\,.$$
\endproclaim

We refer to \cite{8} for further discussion of `local to global' formulas.

\s\smc

 B. B.: Institute of Mathematics PAN, 

ul \'Sniadeckich 8, 00-950 Warszawa, Poland

(Instytut Matematyczny PAN)

\tt bojarski\@impan.gov.pl 

\s\smc

A. W.: Institute of Mathematics, Warsaw University, 

ul.Banacha 2, 02-097, Warszawa, Poland 

(Instytut Matematyki, Uniwersytet Warszawski)

\tt aweber\@mimuw.edu.pl

\enddocument
\end